\documentclass[12pt]{article}%
\usepackage{amsmath,amssymb}
\usepackage[applemac]{inputenc}
\usepackage{amsmath,amssymb,fullpage}
\usepackage{color}
\usepackage{amsmath}
\usepackage{amsfonts}
\usepackage{amssymb}
\usepackage{graphicx}%
\providecommand{\U}[1]{\protect\rule{.1in}{.1in}}

\newtheorem{definition}{Definition}[section]

\newtheorem{theorem}[definition]{Theorem}

\newtheorem{remark}[definition]{ \it Remark}

\newtheorem{lemma}[definition]{Lemma}

\numberwithin{equation}{section}

\def\1B{\text{1\!\!I}}

\begin{document}

\date{23 December 2018}
\title{Dynamic risk measure for BSVIE with jumps and semimartingale issues}
\author{Nacira Agram$^{1,2}$ }
\maketitle

\begin{abstract}
Risk measure is a fundamental concept in finance and in the insurance
industry, it is used to adjust life insurance rates. In this current paper, we
will study dynamic risk measures by means of backward stochastic Volterra
integral equations (BSVIEs) with jumps. We prove a comparison theorem for such
a type of equations. Since the solution of a BSVIEs is not a semimartingale in
general, we will discuss some particular semimartingale issues.

\end{abstract}

\footnotetext[1]{Department of Mathematics, University of Oslo, P.O. Box 1053
Blindern, N--0316 Oslo, Norway and University of Biskra, Algeria. Email:
\texttt{naciraa@math.uio.no}.\newline}

\footnotetext[2]{This research was carried out with support of the Norwegian
Research Council, within the research project Challenges in Stochastic
Control, Information and Applications (STOCONINF), project number 250768/F20.
\par
Part of this work has been done while the author is visiting the University of
Alberta, Canada and I would like to thank Prof. Hu for the hospitality and his
comments during my talk.}

\paragraph{MSC(2010):\newline}

60H07, 60H20, 60H30, 45D05, 45R05.



\section{Introduction}

Consider the solution couple $(Y,Z)$ of the following nonlinear backward
stochastic differential equation (BSDE)%
\begin{equation}
dY(t)=-F(t,Y(t),Z(t))dt+Z(t)dB(t),\text{ }Y(T)=-\xi. \label{YI}%
\end{equation}
It is known that under mild conditions this equation has a unique solution.
This was first proved by Pardoux and Peng \cite{pp}, for a Lipschitz driver
$F$ and a square integrable random variable $\xi$. Such a type of equation
first appeared in its linear form as an adjoint equation when Bismut \cite{b}
studied the stochastic version of the Pontryagin's maximum principle. Due to
their significant applications to finance and insurance, BSDEs have gained a
lot of interest since $1990$. See e.g. the seminal work by El Karoui
\textit{et al} \cite{epq} for more details. Subsequently, it was discovered
that BSDEs could also be used to represent risk measures. We recall this
connection briefly. First, let us recall the definition of a convex risk
measure, see for instance Föllmer and Schied \cite{fs}, Frittelli and Rosazza
Gianin \cite{fr}.

\begin{definition}
A convex risk measure is a map $\rho:L^{p}\left(  \mathcal{F}_{T}\right)
\rightarrow%
\mathbb{R}
,$ $p\in\lbrack2,\infty]$ that satisfies the following properties:
\end{definition}

\begin{itemize}
\item (Convexity) $\rho(\lambda\varphi_{1}+(1-\lambda)\varphi_{2})\leq
\lambda\rho(\varphi_{1})+(1-\lambda)\rho(\varphi_{2})$ for all $\lambda
\in\lbrack0,1]$ and all $\varphi_{1},\varphi_{2}\in L^{p}\left(
\mathcal{F}_{T}\right)  $.

\item (Monotonicity) If $\varphi_{1}\leq\varphi_{2},$ then $\rho(\varphi
_{1})\geq\rho(\varphi_{2})$.

\item (Translation invariance) $\rho(\varphi+a)=\rho(\varphi)-a$ for all
$\varphi\in L^{p}\left(  \mathcal{F}_{T}\right)  $ and all constants $a$.

\item $\rho(0)=0$.
\end{itemize}

The last axiom is known as a normalisation and usually it is assumed for
convenience.\newline The construction of risk measures from solutions to BSDEs
is given as follows: \newline Let $Y(t)=Y^{-\xi}(t)$ be the first component of
the solution of the BSDE (\ref{YI}). Assume that the driver $F(t,y,z)$ of the
BSDE (\ref{YI}) does not depend on $y$ and that $z\mapsto F(t,z)$ is convex
for all $t$. Then
\[
\rho(\xi):=Y^{-\xi}(0)
\]
defines a convex risk measure. This shows how crucial is the choice of the
functional $F$. This relation extends to the concept of \emph{dynamic} risk
measure. See for example Barrieu and El Karoui \cite{be}, Quenez and Sulem
\cite{qs1}. Peng \cite{Peng}, Rosazza Gianin \cite{g} have used the notion of
the nonlinear expectation ($g$-BSDE) to represent a dynamic risk
measures.\newline The connection above has also been generalised to BSDEs with
jumps by Quenez and Sulem \cite{qs} and Øksendal and Sulem \cite{OS3}. Yong
\cite{Yong_risk} extended the notion of dynamic risk measures by means of
BSDEs to a specific class of dynamic risk measures that arise as a solutions
to the so-called \emph{BSVIEs}, as follows: Consider a position (wealth)
process $\psi\in L_{\mathcal{F}_{T}}^{2}(0,T)$ instead of a random variable
$\xi$ and let the couple $(Y(t),Z(t,s))$ be the solution of the following
BSVIE
\begin{equation}
Y(t)=-\psi\left(  t\right)  +%
{\textstyle\int_{t}^{T}}
g(t,s,Y(s),Z(t,s))ds-%
{\textstyle\int_{t}^{T}}
Z(t,s)dB(s),t\in\lbrack0,T]. \label{y1}%
\end{equation}
More generally, he considers the BSVIE%
\begin{equation}
Y(t)=-\psi\left(  t\right)  +%
{\textstyle\int_{t}^{T}}
g(t,s,Y(s),Z(t,s),Z(s,t))ds-%
{\textstyle\int_{t}^{T}}
Z(t,s)dB(s),t\in\left[  0,T\right]  . \label{y2}%
\end{equation}
We may remark here, the difference between the two above equations is that the
driver in (\ref{y2}) depends on both $Z(t,s)$ and $Z(s,t).$ This type of
equations appeared in its linear form when Yong \cite{Yong1} proved the
duality principle between linear forward stochastic Volterra integral
equations (SVIEs) and linear BSVIEs.\newline Moreover, differentiability of
BSDEs and BSVIEs were used for the purpose of capital allocation for risk
measures by Kromer and Overbeck \cite{KO}. We recall now the definition of the
dynamic risk measure by means of BSVIEs from Yong \cite{Yong_risk}.

\begin{definition}
A map $\rho:[0,T]\times L_{\mathcal{F}_{T}}^{2}(0,T)\rightarrow L_{\mathbb{F}%
}^{2}(0,T)$ is called a dynamic risk measure if the following hold:
\end{definition}

\begin{itemize}
\item (Past independence) For any $\psi\left(  \cdot\right)  ,\psi^{\prime
}\left(  \cdot\right)  \in L_{\mathcal{F}_{T}}^{2}(0,T),$ if%
\[
\psi\left(  s\right)  =\psi^{\prime}\left(  s\right)  ,\text{ \ \ a.s.,
}\omega\in\Omega,\text{ }s\in\lbrack t,T],
\]
for some $t\in\lbrack0,T)$, then%
\[
\rho(t;\psi\left(  \cdot\right)  )=\rho(t;\psi^{\prime}\left(  \cdot\right)
),\text{ \ \ a.s., }\omega\in\Omega.
\]

\item (Monotonicity) For any $\psi\left(  \cdot\right)  ,\psi^{\prime}\left(
\cdot\right)  \in L_{\mathcal{F}_{T}}^{2}(0,T),$ if%
\[
\psi\left(  s\right)  \leq\psi^{\prime}\left(  s\right)  ,\text{ \ \ a.s.,
}\omega\in\Omega,\text{ }s\in\lbrack t,T],
\]
for some $t\in\lbrack0,T)$, then%
\[
\rho(t;\psi\left(  \cdot\right)  )\geq\rho(t;\psi^{\prime}\left(
\cdot\right)  ),\text{ \ \ a.s., }\omega\in\Omega.
\]

\end{itemize}

\begin{definition}
A dynamic risk measure $\rho:L_{\mathcal{F}_{T}}^{2}(0,T)\rightarrow
L_{\mathbb{F}}^{2}(0,T)$ is called a convex risk measure if
\end{definition}

\begin{itemize}
\item (Convexity) For any $\psi\left(  \cdot\right)  ,\psi^{\prime}\left(
\cdot\right)  \in L_{\mathcal{F}_{T}}^{2}(0,T)$ and $\lambda\in\lbrack0,1],$%
\[
\rho(t;\lambda\psi\left(  \cdot\right)  +(1-\lambda)\psi^{\prime}\left(
\cdot\right)  )\leq\lambda\rho(t;\psi\left(  \cdot\right)  )+(1-\lambda
)\rho(t;\psi^{\prime}\left(  \cdot\right)  ),\text{ \ \ a.s., }\omega\in
\Omega,t\in\lbrack0,T].
\]

\item (Translation invariance) For any $\psi\left(  \cdot\right)  $ and any
constant $a$ it holds that
\[
\rho(t;\psi\left(  \cdot\right)  +a)=\rho(t;\psi\left(  \cdot\right)
)-a,\text{ \ \ a.s., }\omega\in\Omega,\text{ }s\in\lbrack t,T],
\]

\end{itemize}

This type of equation can also be used as a model for recursive utility,
because the equivalent formulation of (\ref{y1}) is%

\[%
\begin{array}
[c]{cc}%
Y(t) & =\mathbb{E[-}\psi(t)+%
{\textstyle\int_{t}^{T}}
g(t,s,Y(s),Z(t,s))ds|\mathcal{F}_{t}].
\end{array}
\]
This can be regarded as an extension of the classical recursive utility
concept of Duffie and Epstein \cite{DE} to systems with memory. For example,
for the consumption process $c(t)\geq0$, we consider its recursive utility
process $Y(t)$ defined by
\[%
\begin{array}
[c]{cc}%
Y(t) & =\mathbb{E[-}\psi(t)+%
{\textstyle\int_{t}^{T}}
g(t,s,Y(s),Z(t,s),c(s))ds|\mathcal{F}_{t}].
\end{array}
\]
For more details about BSVIEs and their applications, we refer to Yong
\textit{et al} \cite{Yong1}, \cite{Yong 2}, \cite{Yong}\ and to Agram
\textit{et al} \cite{AOY},\cite{AOY2}, \cite{AO}. In Wang and Yong \cite{wy}
comparison theorems for different classes of BSVIEs have been proved. In all
the above works, only BSVIEs driven by Brownian motion are considered.

In the current paper we are interested in BSVIEs with jumps of the form
\begin{equation}%
\begin{array}
[c]{c}%
Y(t)=-\psi\left(  t\right)  +\int_{t}^{T}g(t,s,Y(s),Z(t,s),K(t,s,\cdot
))ds-\int_{t}^{T}Z(t,s)dB(s)\\
-\int_{t}^{T}\int_{%
\mathbb{R}
_{0}}K(t,s,\zeta)\tilde{N}(ds,d\zeta),t\in\left[  0,T\right]  ,
\end{array}
\label{yi}%
\end{equation}
where $B$ is a standard Brownian motion and $\tilde{N}$ is an independent
compensated Poisson random measure (see below for definitions). For a
Lipschitz driver, Agram \textit{et al} \cite{AOY} have proved that there
exists a unique solution of equation (\ref{yi}). Under weaker assumptions
(non-Lipschitz driver), we refer to Wang and Zhang \cite{WZ} and to Ren
\cite{ren}.\newline The general nature of the BSVIEs does not allow us to
write explicitly the solution of a linear BSVIEs in general. However, Hu and
Øksendal \cite{HO} have obtained a closed solution formula for a special class
of linear BSVIEs with jumps. Using their result, we will prove comparison
theorems for BSVIEs with jumps.\newline

The main contributions in this paper is the extension to jumps of the
comparison theorems in Wang and Yong \cite{wy} and also the dynamic measures
by BSVIEs with jumps of the paper by Yong \cite{Yong_risk}.

Moreover, due to the dependence on $t$ in the coefficients, such BSVIEs are
complicated to deal with. It is not even clear if the solution of a BSVIE is a
semimartingale in general. We will discuss some particular cases where the
solution can be a semimartingale.\newline\newline Here is an outline of our
paper:\newline In Section 2 we give some preliminaries on BSVIEs with jumps.
Then in Section 3 we prove a comparison theorem for BSVIEs with jumps. In
Section 4 we study dynamic risk measures by means of BSVIEs with jumps.
Finally, in Section 5 we discuss some semimartingale issues for BSVIEs.

\section{Preliminaries on BSVIEs with jumps}

Let $(%
\Omega
,\mathcal{F},\mathbb{P})$ be a given probability space with filtration
$\mathbb{F}=(\mathcal{F}_{t})_{t\geq0}$ generated by a one-dimensional
Brownian motion $B$ and an independent Poisson random measure $N(dt,d\zeta)$.
Let $\nu(d\zeta)dt$ denote the Lévy measure of $N$, and let $\tilde
{N}(dt,d\zeta)$ denote the compensated Poisson random measure $N(dt,d\zeta
)-\nu(d\zeta)dt$.\newline Put $\triangle:=\{\left(  t,s\right)  \in\left[
0,T\right]  ^{2}:t\leq s\}$. We define the following sets:

\begin{itemize}
\item $L_{y}^{2}$ consists of the $\mathbb{F}$-adapted càdlàg processes
$Y:[0,T]\times\Omega\rightarrow\mathbb{R}$ equipped with the norm
\[
\parallel Y\parallel_{L_{y}^{2}}^{2}:=\mathbb{E[}%
{\textstyle\int_{0}^{T}}
|Y(t)|^{2}dt]<\infty.
\]

\item $L_{z}^{2}$ consists of the $\mathbb{F}$-predictable processes
\[
Z:\triangle\times\Omega\rightarrow\mathbb{R},
\]
such that $\mathbb{E[}%
{\textstyle\int_{0}^{T}}
{\textstyle\int_{t}^{T}}
\left\vert Z(t,s)\right\vert ^{2}dsdt]<\infty$ with $s\mapsto Z(t,s)$ being
$\mathbb{F}$-predictable on $[t,T].$ We equip $L_{z}^{2}$ with the norm
\[
\parallel Z\parallel_{L_{z}^{2}}^{2}:=\mathbb{E[}%
{\textstyle\int_{0}^{T}}
{\textstyle\int_{t}^{T}}
\left\vert Z(t,s)\right\vert ^{2}dsdt].
\]
{}

\item $L_{\nu}^{2}$ consists of all Borel functions $K:%
\mathbb{R}
_{0}\rightarrow\mathbb{R},$ such that
\[
\parallel K\parallel_{L_{\nu}^{2}}^{2}:=%
{\textstyle\int_{\mathbb{R}_{0}}}
K(t,s,\zeta)^{2}\nu(d\zeta)<\infty.
\]

\item $H_{\nu}^{2}$ consists of $\mathbb{F}$-predictable
processes$\ K:\triangle\times%
\mathbb{R}
_{0}\times\Omega\rightarrow\mathbb{R},$ such that
\[
\mathbb{E[}%
{\textstyle\int_{0}^{T}}
{\textstyle\int_{t}^{T}}
{\textstyle\int_{\mathbb{R}_{0}}}
|K(t,s,\zeta)|^{2}\nu(d\zeta)dsdt]<\infty
\]
and $s\mapsto K\left(  t,s,\cdot\right)  $ being $\mathbb{F}$-predictable on
$[t,T].$ We equip $H_{\nu}^{2}$ with the norm
\[
\parallel K\parallel_{H_{\nu}^{2}}^{2}:=\mathbb{E[}%
{\textstyle\int_{0}^{T}}
{\textstyle\int_{t}^{T}}
{\textstyle\int_{\mathbb{R}_{0}}}
|K(t,s,\zeta)|^{2}\nu(d\zeta)dsdt].
\]

\item Let $L_{\mathcal{F}_{T}}^{2}[0,T]$ be the space of all processes
$\psi:[0,T]\times\Omega\rightarrow\mathbb{R}$ and $\psi$ is $\mathcal{F}_{T}%
$-measurable for all $t\in\lbrack0,T],$ such that%
\[
||\psi||_{L_{\mathcal{F}_{T}}^{2}[0,T]}^{2}=\mathbb{E[}%
{\textstyle\int_{0}^{T}}
|\psi(t)|^{2}dt]<\infty.
\]

\item $L_{\mathbb{F}}^{2}[0,T]$ is the space of all $\psi\in L_{\mathcal{F}%
_{T}}^{2}[0,T]$ that are $\mathbb{F}$-adapted.
\end{itemize}

\bigskip

\noindent We are interested in the BSVIE $(Y,Z,K)\in L_{y}^{2}\times L_{z}%
^{2}\times H_{\nu}^{2}$, given by
\begin{equation}%
\begin{array}
[c]{c}%
Y(t)=\psi\left(  t\right)  +\int_{t}^{T}g(t,s,Y(s),Z(t,s),K(t,s,\cdot
))ds-\int_{t}^{T}Z(t,s)dB(s)\\
-\int_{t}^{T}\int_{%
\mathbb{R}
_{0}}K(t,s,\zeta)\tilde{N}(ds,d\zeta),t\in\left[  0,T\right]  ,
\end{array}
\label{a3}%
\end{equation}
where the following conditions are satisfied:

\begin{itemize}
\item[(i)] The driver $g:\triangle\times%
\mathbb{R}
^{2}\times L_{\nu}^{2}\times\Omega\rightarrow%
\mathbb{R}
$ satisfies the following integrability condition:%

\begin{equation}
\mathbb{E[}%
{\textstyle\int\nolimits_{0}^{T}}
(%
{\textstyle\int\nolimits_{t}^{T}}
g(t,s,0,0,0)ds)^{2}dt]<+\infty. \label{integr}%
\end{equation}

\item[(ii)] The driver $g$ satisfies the following Lipschitz
condition:\newline There exists a constant $C>0$, such that, for all
$(t,s)\in\triangle,$%
\begin{equation}%
\begin{array}
[c]{l}%
\left\vert g(t,s,y,z,k(\cdot))-g(t,s,y^{\prime},z^{\prime},k^{\prime}%
(\cdot))\right\vert \\
\leq C(\left\vert y-y^{\prime}\right\vert +\left\vert z-z^{\prime}\right\vert
+(\int_{%
\mathbb{R}
_{0}}\left\vert k(\zeta)-k^{\prime}\left(  \zeta\right)  \right\vert ^{2}%
\nu(d\zeta))^{\frac{1}{2}}),
\end{array}
\label{lipsc}%
\end{equation}
for all $y,y^{\prime},z,z^{\prime}\in%
\mathbb{R}
,k(\cdot),k^{\prime}(\cdot)\in L_{\nu}^{2}.$

\item[(iii)] The terminal condition $\psi\left(  \cdot\right)  \in
L_{\mathcal{F}_{T}}^{2}[0,T].$
\end{itemize}

\noindent For this equation the following is known:

\begin{theorem}
[Agram \textit{et al} \cite{AOY}]\label{exi-uni} Under the above assumptions,
there exists a unique solution $(Y,Z,K)\in L_{y}^{2}\times L_{z}^{2}\times
H_{\nu}^{2}$ of the BSVIEs with jumps $(\ref{a3}),$ with
\[
\left\Vert \left(  Y,Z,K\right)  \right\Vert _{L_{y}^{2}\times L_{z}^{2}\times
H_{\nu}^{2}}^{2}\leq C\mathbb{E}[\left\vert \psi\left(  t\right)  \right\vert
^{2}+(%
{\textstyle\int_{t}^{T}}
g\left(  t,s,0,0,0\right)  ds)^{2}].
\]

\end{theorem}

\noindent In the next section we will prove a comparison theorem for some
BSVIEs by using the linearisation of the solutions.

\section{Comparison theorem for BSVIEs with jumps}

In order to be able to prove the comparison theorem for BSVIEs with jumps, let
us first recall the closed formula for linear BSVIEs with jumps. Consider the
solution triplet $(Y,Z,K)\in L_{y}^{2}\times L_{z}^{2}\times H_{\nu}^{2}$ of
the following linear BSVIE%
\begin{equation}%
\begin{array}
[c]{c}%
Y(t):=\psi\left(  t\right)  +%
{\textstyle\int_{t}^{T}}
[\alpha(t,s)Y(s)+\beta\left(  s\right)  Z(t,s)+%
{\textstyle\int_{\mathbb{R}_{0}}}
\theta\left(  s,\zeta\right)  K(t,s,\zeta)\nu\left(  d\zeta\right)  ]ds\\
-%
{\textstyle\int_{t}^{T}}
Z(t,s)dB(s)-%
{\textstyle\int_{t}^{T}}
{\textstyle\int_{\mathbb{R}_{0}}}
K(t,s,\zeta)\tilde{N}(ds,d\zeta),t\in\left[  0,T\right]  ,
\end{array}
\label{lineary}%
\end{equation}
where $\psi\left(  \cdot\right)  \in L_{\mathcal{F}_{T}}^{2}[0,T],$
$\alpha(t,s);0\leq t\leq s\leq T\,$ and $(\beta(s),\theta(s,\zeta);0\leq s\leq
T\,,\zeta\in{\mathbb{\ R}}_{0}$ are given (deterministic) measurable functions
of $t$, $s$ and $\zeta$, with values in $\mathbb{R}$. For simplicity we assume
that these functions are bounded, and we assume that there exists
$\varepsilon>0$ such that $\theta(s,\zeta)\geq-1+\varepsilon$ for all
$s,\zeta$.\newline Define the measure $\mathbb{Q}$ by
\[
d\mathbb{Q}=M(T)d\mathbb{P}\text{ on }\mathcal{F}_{T},
\]
where
\begin{align*}
M(t)  &  :=\exp(%
{\textstyle\int_{0}^{t}}
\beta(s)dB(s)-\tfrac{1}{2}%
{\textstyle\int_{0}^{t}}
\beta^{2}(s)ds\\
&  +%
{\textstyle\int_{0}^{t}}
{\textstyle\int_{\mathbb{R}_{0}}}
\ln(1+\theta(s,\zeta))\tilde{N}(ds,d\zeta)\\
&  +%
{\textstyle\int_{0}^{t}}
{\textstyle\int_{\mathbb{R}_{0}}}
\{\ln(1+\theta(s,\zeta))-\theta(s,\zeta)\}\nu(d\zeta)ds).
\end{align*}
Then under the new measure $\mathbb{Q}$ the process
\[
B_{\mathbb{Q}}(t):=B(t)-%
{\textstyle\int_{0}^{t}}
\beta(s)ds\,,\quad0\leq t\leq T\,.
\]
is a Brownian motion, and the random measure
\[
\tilde{N}_{\mathbb{Q}}(dt,d\zeta):=\tilde{N}(dt,d\zeta)-\theta(t,\zeta
)\nu(d\zeta)dt
\]
is the $\mathbb{Q}$-compensated Poisson random measure of $N(\cdot,\cdot)$, in
the sense that the process%

\[
\tilde{N}_{\gamma}(t):=%
{\textstyle\int_{0}^{t}}
{\textstyle\int_{\mathbb{R}_{0}}}
\gamma(s,\zeta)\tilde{N}_{\mathbb{Q}}(ds,d\zeta)
\]
is a local $\mathbb{Q}$-martingale, for all predictable processes
$\gamma(t,\zeta)$ such that
\[%
{\textstyle\int_{0}^{T}}
{\textstyle\int_{\mathbb{R}_{0}}}
\gamma^{2}(t,\zeta)\theta^{2}(t,\zeta)\nu(d\zeta)dt<\infty.
\]
By this change of measure the equation (\ref{lineary}) is equivalent to%

\[
Y(t)=\psi(t)+%
{\textstyle\int_{t}^{T}}
\alpha(t,s)Y(s)ds-%
{\textstyle\int_{t}^{T}}
Z(t,s)dB_{\mathbb{Q}}(s)-%
{\textstyle\int_{t}^{T}}
{\textstyle\int_{\mathbb{R}_{0}}}
K(t,s,\zeta)\tilde{N}_{\mathbb{Q}}(ds,d\zeta).
\]
For all $0\leq t\leq r\leq T,$ define
\[
\alpha^{(1)}(t,r)=\alpha(t,r)\,,\quad\alpha^{(2)}(t,r)=%
{\textstyle\int_{t}^{r}}
\alpha(t,s)\alpha(s,r)ds
\]
and inductively
\[
\alpha^{(n)}(t,r)=%
{\textstyle\int_{t}^{r}}
\alpha^{(n-1)}(t,s)\alpha(s,r)ds\,,n=3,4,\cdots\,.
\]
Note that if $|\alpha(t,r)|\leq C$ (constant) for all $t,r$, then by induction
on $n\in%
\mathbb{N}
$
\[
|\alpha^{(n)}(t,r)|\leq\tfrac{C^{n}T^{n}}{n!},
\]
for all $t,r,n$. Put,
\[
\Psi(t,r):=\Sigma_{n=1}^{\infty}|\alpha^{(n)}(t,r)|<\infty,
\]
\bigskip for all $t,r$. Then we have the following:

\begin{theorem}
[Hu and Øksendal \cite{HO}]\label{closed-form} The first component $Y(t)$ of
the solution of the linear BSVIE with jumps (\ref{lineary}) is given by
\begin{equation}
Y(t)=\frac{\mathbb{E}\left[  \frac{d\mathbb{Q}}{d\mathbb{P}}\{\psi(t)+%
{\textstyle\int_{t}^{T}}
\Phi(t,r)\psi(r)dr\}|\mathcal{F}_{t}\right]  }{\mathbb{E}\left[
\frac{d\mathbb{Q}}{d\mathbb{P}}|\mathcal{F}_{t}\right]  }.\nonumber
\end{equation}

\end{theorem}

\begin{remark}
In the previous theorem, we can get the same result by considering predictable
processes $(\beta(s),\theta(s,\zeta);0\leq s\leq T\,,\zeta\in{\mathbb{\ R}%
}_{0})$ instead of deterministic functions since the Girsanov change of
measure theorem is still valid.
\end{remark}

\noindent We now state and prove the comparison theorem.

\begin{theorem}
[Comparison Theorem]\label{comp 2} For $i=1,2$, let $g_{i}:\triangle
\times\mathcal{%
\mathbb{R}
}^{2}\times L_{\nu}^{2}\times\Omega\mathcal{\rightarrow%
\mathbb{R}
}$ and $\psi_{1}(t),\psi_{2}(t)\in L_{\mathcal{F}_{T}}^{2}[0,T]$ and let
$(Y^{i},Z^{i},K^{i})_{i=1,2}$ be the solutions of%
\[%
\begin{array}
[c]{c}%
Y^{i}(t)=\psi^{i}\left(  t\right)  +%
{\textstyle\int_{t}^{T}}
g_{i}\left(  t,s,Y^{i}\left(  s\right)  ,Z^{i}\left(  t,s\right)
,K^{i}\left(  t,s,\cdot\right)  \right)  ds-%
{\textstyle\int_{t}^{T}}
Z^{i}(t,s)dB(s)\\
-%
{\textstyle\int_{t}^{T}}
{\textstyle\int_{\mathbb{R}_{0}}}
K^{i}(t,s,\zeta)\tilde{N}(ds,d\zeta),t\in\left[  0,T\right]  .
\end{array}
\]
Assume that the drivers $(g_{i})_{i=1,2}$ are Lipschitz and satisfy
\begin{equation}
g_{1}\left(  t,s,y^{2},z^{2},k^{2}\right)  \geq g_{2}\left(  t,s,y^{2}%
,z^{2},k^{2}\right)  ,\forall t,\mathbb{P}\text{-a.s.,} \label{g}%
\end{equation}
and that there exists a bounded predictable process $\theta\left(
s,t,\zeta\right)  $ and a $\Pi(\cdot)\in L_{\nu}^{2}$ such that $ds\otimes
d\mathbb{P\otimes}\nu(d\zeta)$-a.s.,
\begin{equation}
\theta\left(  s,\zeta\right)  \geq-1+\varepsilon\text{ and }|\theta\left(
s,\zeta\right)  |\leq\Pi\left(  \zeta\right)  , \label{gcomp}%
\end{equation}

and the following inequality holds
\begin{equation}%
\begin{array}
[c]{c}%
g_{1}\left(  t,s,Y^{2}\left(  s\right)  ,Z^{2}\left(  t,s\right)
,K^{1}\left(  t,s,\cdot\right)  \right)  -g_{1}\left(  t,s,Y^{2}\left(
s\right)  ,Z^{2}\left(  t,s\right)  ,K^{2}\left(  t,s,\cdot\right)  \right) \\
\geq%
{\textstyle\int_{\mathbb{R}_{0}}}
\theta\left(  s,\zeta\right)  (K^{1}(t,s,\zeta)-K^{2}(t,s,\zeta))\nu(d\zeta).
\end{array}
\label{*}%
\end{equation}

Moreover, assume that the driver $g$ is increasing on $y$, such that%
\begin{equation}
g\left(  t,s,y^{1},z,k\right)  \geq g\left(  t,s,y^{2},z,k\right)  ,\text{ if
}y^{1}\geq y^{2},\text{\ }\forall t\text{,}\mathbb{P}\text{-a.s.,}
\label{incy}%
\end{equation}
and%
\begin{equation}
\psi_{1}(t)\geq\psi_{2}(t)\text{ for each }t\in\left[  0,T\right]  ,\text{
}\mathbb{P}\text{-a.s.} \label{xi}%
\end{equation}

Then $Y^{1}(t)\geq Y^{2}(t)$ $\mathbb{P}$-a.s. for each $t$.\newline
\end{theorem}

\noindent{Proof}\quad We set%
\begin{align*}
\hat{\psi}  &  =\psi_{1}-\psi_{2},\\
\hat{Y}  &  =Y^{1}-Y^{2},\text{ }\hat{Z}=Z^{1}-Z^{2},\text{ }\hat{K}%
=K^{1}-K^{2},
\end{align*}
we have%
\[%
\begin{array}
[c]{c}%
\hat{Y}(t)=\hat{\psi}\left(  t\right)  +\int_{t}^{T}[g_{1}\left(
t,s,Y^{1}\left(  s\right)  ,Z^{1}\left(  t,s\right)  ,K^{1}\left(
t,s,\cdot\right)  \right) \\
\quad\quad\quad\quad\quad-g_{2}\left(  t,s,Y^{2}\left(  s\right)
,Z^{2}\left(  t,s\right)  ,K^{2}\left(  t,s,\cdot\right)  \right)  ]ds\\
-\int_{t}^{T}\hat{Z}(t,s)dB(s)-\int_{t}^{T}\int_{\mathbb{R}_{0}}\hat
{K}(t,s,\zeta)\tilde{N}(ds,d\zeta),t\in\left[  0,T\right]  .
\end{array}
\]
Note that%
\begin{align*}
&  g_{1}\left(  t,s,Y^{1}\left(  s\right)  ,Z^{1}\left(  t,s\right)
,K^{1}\left(  t,s,\cdot\right)  \right)  -g_{2}\left(  t,s,Y^{2}\left(
s\right)  ,Z^{2}\left(  t,s\right)  ,K^{2}\left(  t,s,\cdot\right)  \right) \\
&  =g_{1}\left(  t,s,Y^{1}\left(  s\right)  ,Z^{1}\left(  t,s\right)
,K^{1}\left(  t,s,\cdot\right)  \right)  -g_{1}\left(  t,s,Y^{2}\left(
s\right)  ,Z^{1}\left(  t,s\right)  ,K^{1}\left(  t,s,\cdot\right)  \right) \\
&  +g_{1}\left(  t,s,Y^{2}\left(  s\right)  ,Z^{1}\left(  t,s\right)
,K^{1}\left(  t,s,\cdot\right)  \right)  -g_{1}\left(  t,s,Y^{2}\left(
s\right)  ,Z^{2}\left(  t,s\right)  ,K^{1}\left(  t,s,\cdot\right)  \right) \\
&  +g_{1}\left(  t,s,Y^{2}\left(  s\right)  ,Z^{2}\left(  t,s\right)
,K^{1}\left(  t,s,\cdot\right)  \right)  -g_{1}\left(  t,s,Y^{2}\left(
s\right)  ,Z^{2}\left(  t,s\right)  ,K^{2}\left(  t,s,\cdot\right)  \right) \\
&  g_{1}\left(  t,s,Y^{2}\left(  s\right)  ,Z^{2}\left(  t,s\right)
,K^{2}\left(  t,s,\cdot\right)  \right)  -g_{2}\left(  t,s,Y^{2}\left(
s\right)  ,Z^{2}\left(  t,s\right)  ,K^{2}\left(  t,s,\cdot\right)  \right) \\
&  \geq\alpha\left(  t,s\right)  \hat{Y}\left(  s\right)  +\beta\left(
s\right)  \hat{Z}(t,s)+%
{\textstyle\int_{\mathbb{R}_{0}}}
\theta\left(  s,\zeta\right)  \hat{K}(t,s,\zeta)\nu\left(  d\zeta\right)  ,
\end{align*}
where%
\[
\alpha\left(  t,s\right)  =\tfrac{g_{1}\left(  t,s,Y^{1}\left(  s\right)
,Z^{1}\left(  t,s\right)  ,K^{1}\left(  t,s,\cdot\right)  \right)
-g_{1}\left(  t,s,Y^{2}\left(  s\right)  ,Z^{1}\left(  t,s\right)
,K^{1}\left(  t,s,\cdot\right)  \right)  }{\hat{Y}\left(  s\right)
}\mathbf{1}_{\left\{  \hat{Y}\left(  s\right)  \neq0\right\}  },
\]
and%
\[
\beta\left(  s\right)  =\tfrac{g_{1}(t,s,Y^{2}\left(  s\right)  ,Z^{1}\left(
t,s\right)  ,K^{1}\left(  t,s,\cdot\right)  )-g_{1}(t,s,Y^{2}\left(  s\right)
,Z^{2}\left(  t,s\right)  ,K^{1}\left(  t,s,\cdot\right)  )}{\hat{Z}\left(
t,s\right)  }\mathbf{1}_{\{\hat{Z}\left(  t,s\right)  \neq0\}}.
\]
Hence%
\begin{equation}%
\begin{array}
[c]{c}%
\hat{Y}(t)\geq\hat{\psi}\left(  t\right)  +%
{\textstyle\int_{t}^{T}}
[\alpha\left(  t,s\right)  \hat{Y}\left(  s\right)  +\beta\left(  s\right)
\hat{Z}(t,s)+%
{\textstyle\int_{\mathbb{R}_{0}}}
\theta\left(  s,t,\zeta\right)  \hat{K}(t,s,\zeta)\nu\left(  d\zeta\right)
]ds\\
-%
{\textstyle\int_{t}^{T}}
\hat{Z}(t,s)dB(s)-%
{\textstyle\int_{t}^{T}}
{\textstyle\int_{\mathbb{R}_{0}}}
\hat{K}(t,s,\zeta)\tilde{N}(ds,d\zeta),t\in\left[  0,T\right]  .
\end{array}
\label{nbsvie}%
\end{equation}
We proceed as in Theorem \ref{closed-form} and obtain
\begin{equation}
\hat{Y}(t)\geq\frac{\mathbb{E}\left[  \frac{d\mathbb{Q}}{d\mathbb{P}}%
\{\hat{\psi}(t)+%
{\textstyle\int_{t}^{T}}
\Phi(t,r)\hat{\psi}(r)dr\}|\mathcal{F}_{t}\right]  }{\mathbb{E}\left[
\frac{d\mathbb{Q}}{d\mathbb{P}}|\mathcal{F}_{t}\right]  },\nonumber
\end{equation}
where $\Phi(t,r)=\sum_{n=1}^{\infty}\alpha^{(n)}(t,r).$

\noindent By hypothesis (\ref{incy}), we get that $\alpha(t,r)\geq0$ for all
$t,r$. Hence $\Phi(t,r)\geq0$ for all $t,r$ together with (\ref{xi}) implies
that $\hat{Y}(t)\geq0$ $\mathbb{P}$-a.s.\qquad$\qquad\qquad\qquad\qquad
\qquad\qquad\qquad\qquad\qquad\qquad\qquad\qquad\qquad\square$\newline\newline
The following particular case is essential for the next section.

\subsection{ BSVIEs with drivers independent of $Y$}

In this subsection, we are interested in BSVIEs with drivers $g$ independent
of $Y$, as follows:
\begin{equation}%
\begin{array}
[c]{c}%
Y(t)=-\psi\left(  t\right)  +\int_{t}^{T}g(t,s,Z(t,s),K(t,s,\cdot))ds-\int
_{t}^{T}Z(t,s)dB(s)\\
-\int_{t}^{T}\int_{%
\mathbb{R}
_{0}}K(t,s,\zeta)\tilde{N}(ds,d\zeta),t\in\left[  0,T\right]  .
\end{array}
\label{-g}%
\end{equation}
We impose the following set of assumptions.

\begin{itemize}
\item The driver $g:\triangle\times%
\mathbb{R}
\times L_{\nu}^{2}\times\Omega\rightarrow%
\mathbb{R}
$ satisfies
\end{itemize}

\begin{equation}
\mathbb{E[}%
{\textstyle\int\nolimits_{0}^{T}}
(%
{\textstyle\int\nolimits_{t}^{T}}
g(t,s,0,0)ds)^{2}dt]<+\infty. \label{in}%
\end{equation}

\begin{itemize}
\item There exists a constant $C>0$, such that, for all $(t,s)\in\triangle,$%
\begin{equation}%
\begin{array}
[c]{l}%
\left\vert g(t,s,z,k(\cdot))-g(t,s,z^{\prime},k^{\prime}(\cdot))\right\vert \\
\leq C(\left\vert z-z^{\prime}\right\vert +(\int_{%
\mathbb{R}
_{0}}\left\vert k(\zeta)-k^{\prime}\left(  \zeta\right)  \right\vert ^{2}%
\nu(d\zeta))^{\frac{1}{2}}),
\end{array}
\label{li}%
\end{equation}
for all $z,z^{\prime}\in%
\mathbb{R}
,k(\cdot),k^{\prime}(\cdot)\in L_{\nu}^{2}.$

\item The terminal value $\psi\left(  \cdot\right)  \in L_{\mathcal{F}_{T}%
}^{2}[0,T].$
\end{itemize}

\noindent We know by Theorem \ref{exi-uni} there exists a unique solution
$(Y,Z,K)\in L_{y}^{2}\times L_{z}^{2}\times H_{\nu}^{2}$ for the BSVIE
(\ref{-g}). Consider the following the particular linear BSVIE with jumps%
\begin{equation}%
\begin{array}
[c]{c}%
Y(t):=-\psi\left(  t\right)  +%
{\textstyle\int_{t}^{T}}
[\beta\left(  s\right)  Z(t,s)+%
{\textstyle\int_{\mathbb{R}_{0}}}
\theta\left(  s,\zeta\right)  K(t,s,\zeta)\nu\left(  d\zeta\right)  ]ds\\
-%
{\textstyle\int_{t}^{T}}
Z(t,s)dB(s)-%
{\textstyle\int_{t}^{T}}
{\textstyle\int_{\mathbb{R}_{0}}}
K(t,s,\zeta)\tilde{N}(ds,d\zeta),t\in\left[  0,T\right]  ,
\end{array}
\label{linearyy}%
\end{equation}
where $\psi\left(  \cdot\right)  \in L_{\mathcal{F}_{T}}^{2}[0,T]$ and
$(\beta(s),\theta(s,\zeta);0\leq s\leq T\,,\zeta\in{\mathbb{\ R}}_{0}$ are
given (deterministic) measurable functions of $s$ and $\zeta$, with values in
$\mathbb{R}$. For simplicity we assume that these functions are bounded, and
we assume that there exists $\varepsilon>0$ such that $\theta(s,\zeta
)\geq-1+\varepsilon$ for all $s,\zeta$.

\noindent Since the following results are just particular cases of the
previous theorems, we will state them without proofs.

\begin{lemma}
\label{cf} The part of the solution $Y(t)$ of the linear BSVIE with jumps
(\ref{linearyy}) can be given on its closed formula, as
\[
Y(t)=\frac{-\mathbb{E}\left[  \frac{d\mathbb{Q}}{d\mathbb{P}}\psi
(t)|\mathcal{F}_{t}\right]  }{\mathbb{E}\left[  \frac{d\mathbb{Q}}%
{d\mathbb{P}}|\mathcal{F}_{t}\right]  }\,.
\]

\end{lemma}

\noindent We can also get the comparison theorem for BSVIE with jumps of type
(\ref{-g}).

\begin{theorem}
[Comparison Theorem ]\label{comp} For $i=1,2$, let $g_{i}:\triangle
\times\mathcal{%
\mathbb{R}
}\times L_{\nu}^{2}\times\Omega\mathcal{\rightarrow%
\mathbb{R}
}$ and $\psi_{1}(t),\psi_{2}(t)\in L_{\mathcal{F}_{T}}^{2}[0,T]$ and let
$(Y^{i},Z^{i},K^{i})_{i=1,2}$ be the solutions of%
\[%
\begin{array}
[c]{c}%
Y^{i}(t)=-\psi^{i}\left(  t\right)  +%
{\textstyle\int_{t}^{T}}
g_{i}\left(  t,s,Z^{i}\left(  t,s\right)  ,K^{i}\left(  t,s,\cdot\right)
\right)  ds-%
{\textstyle\int_{t}^{T}}
Z^{i}(t,s)dB(s)\\
-%
{\textstyle\int_{t}^{T}}
{\textstyle\int_{\mathbb{R}_{0}}}
K^{i}(t,s,\zeta)\tilde{N}(ds,d\zeta),t\in\left[  0,T\right]  .
\end{array}
\]
Assume that the driver $(g_{i})_{i=1,2}$ is Lipschitz and satisfies
\begin{equation}
g_{1}\left(  t,s,z^{2},k^{2}\right)  \geq g_{2}\left(  t,s,z^{2},k^{2}\right)
,\forall t,\mathbb{P}\text{-a.s.,} \label{inqg}%
\end{equation}
and that there exists a bounded predictable process $\theta\left(
s,\zeta\right)  $ and a $\Pi(\cdot)\in L_{\nu}^{2}$ such that $ds\otimes
d\mathbb{P\otimes}\nu(d\zeta)$-a.s.,
\[
\theta\left(  s,\zeta\right)  \geq-1+\varepsilon\text{ and }|\theta\left(
s,\zeta\right)  |\leq\Pi\left(  \zeta\right)  ,
\]

and the following inequality holds
\begin{equation}%
\begin{array}
[c]{c}%
g_{1}\left(  t,s,Z^{2}\left(  t,s\right)  ,K^{1}\left(  t,s,\cdot\right)
\right)  -g_{1}\left(  t,s,Z^{2}\left(  t,s\right)  ,K^{2}\left(
t,s,\cdot\right)  \right) \\
\geq%
{\textstyle\int_{\mathbb{R}_{0}}}
\theta\left(  s,\zeta\right)  (K^{1}(t,s,\zeta)-K^{2}(t,s,\zeta))\nu(d\zeta).
\end{array}
\label{inqjump}%
\end{equation}

Moreover,%
\begin{equation}
\psi_{1}(t)\leq\psi_{2}(t)\text{ for each }t\in\left[  0,T\right]  ,\text{
}\mathbb{P}\text{-a.s.} \label{psi}%
\end{equation}

Then $Y^{1}(t)\geq Y^{2}(t)$ $\mathbb{P}$-a.s. for each $t$.\newline
\end{theorem}

\section{Dynamic risk measure by means of BSVIE}

As we have seen in the introduction, a natural way to construct a dynamic risk
measures by means of a BSVIEs with jumps, is as follows: Define
\[
\rho(t;\psi(\cdot)):=Y^{-\psi(\cdot)}(t),\,\,\,\,\,\text{for all }t\in
\lbrack0,T],
\]
where $Y$ is the first component of the solution $(Y(t),Z(t,s),K(t,s,\zeta))$
of the BSVIE%
\[%
\begin{array}
[c]{c}%
Y(t)=-\psi\left(  t\right)  +\int_{t}^{T}g(t,s,Z(t,s),K(t,s,\cdot))ds-\int
_{t}^{T}Z(t,s)dB(s)\\
-\int_{t}^{T}\int_{%
\mathbb{R}
_{0}}K(t,s,\zeta)\tilde{N}(ds,d\zeta),t\in\left[  0,T\right]  ,
\end{array}
\]
where the terminal condition $\psi\in L_{\mathcal{F}_{T}}^{2}(0,T)$ and the
generator $g:\triangle\times%
\mathbb{R}
\times L_{\nu}^{2}\times\Omega\rightarrow%
\mathbb{R}
$ satisfies the assumptions (\ref{in}), (\ref{li}), (\ref{gcomp}).

\begin{theorem}
If $\rho$ and $g$ are defined as above then $\rho$ is a \emph{convex dynamic
risk measure}, i.e. the following holds:

\begin{description}
\item[(i) Convexity:] Suppose that $(z,k(\cdot))\mapsto g(t,s,z,k(\cdot))$ is
convex, i.e.,%
\begin{align*}
&  g(t,s,\lambda z_{1}+(1-\lambda)z_{2},\lambda k_{1}(\cdot)+(1-\lambda
)k_{2}(\cdot))\\
&  \leq g(t,s,\lambda z_{1},\lambda k_{1}(\cdot))+g(t,s,(1-\lambda
)z_{2},(1-\lambda)k_{2}(\cdot)),
\end{align*}
for all $(t,s)\in\triangle,z_{1},z_{2}\in%
\mathbb{R}
,k_{1}(\cdot),k_{2}(\cdot)\in%
\mathbb{R}
_{0}$ and $\lambda\in\lbrack0,1].$
\end{description}

Then $\psi(\cdot)\mapsto\rho(t;\psi(\cdot))$ is convex, i.e.,%
\[
\rho(t;\lambda\psi_{1}(\cdot)+(1-\lambda)\psi_{2}(\cdot))\leq\lambda
\rho(t;\psi_{1}(\cdot))+(1-\lambda)\rho(t;\psi_{2}(\cdot)),t\in\lbrack0,T].
\]

\begin{description}
\item[(ii) Monotonicity:] If $\psi_{1}(\cdot)\leq\psi_{2}(\cdot)$, then
$\rho(t;\psi_{2}(\cdot))\leq\rho(t;\psi_{1}(\cdot)).$

\item[(iii) Translation invariance:] If $\psi(\cdot)\in L_{\mathcal{F}_{T}%
}^{2}(0,T)$ and a constant $a\in%
\mathbb{R}
$. Then $\rho(t;\psi(\cdot)+a)=\rho(t;\psi(\cdot))-a,\text{ for each }%
t\in\left[  0,T\right]  . $

\item[(iv) Past independence:] If $\psi(\cdot),\psi^{\prime}(\cdot)\in
L_{\mathcal{F}_{T}}^{2}(0,T)$ and $\psi(s)=\psi^{\prime}(s)$ for all
$s\in\lbrack t,T]$ then $\rho(t;\psi(\cdot))=\rho(t,\psi^{\prime}(\cdot)).$
\end{description}
\end{theorem}

\noindent{Proof}\quad(i) Convexity: Fix $\lambda\in(0,1)$ and for all
$\psi_{1}(\cdot),\psi_{2}(\cdot)\in L_{\mathcal{F}_{T}}^{2}(0,T)$. We want to
prove that
\[
\rho(t;\lambda\psi_{1}(\cdot)+(1-\lambda)\psi_{2}(\cdot))\leq\lambda
\rho(t;\psi_{1}(\cdot))+(1-\lambda)\rho(t;\psi_{2}(\cdot)),
\]
i.e.,%

\[
Y^{-(\lambda\psi_{1}(\cdot)+(1-\lambda)\psi_{2}(\cdot))}(t)\leq\lambda
(Y^{-\psi_{1}(\cdot)}(t))+(1-\lambda)(Y^{-\psi_{2}(\cdot)}(t)).
\]
Set $(\hat{Y},\hat{Z},\hat{K})\in L_{y}^{2}\times L_{z}^{2}\times H_{\nu}^{2}$
solution of the following BSVIE with jumps%
\[%
\begin{array}
[c]{c}%
\hat{Y}(t)=-\lambda\psi_{1}(t)-(1-\lambda)\psi_{2}(t)+\int_{t}^{T}%
g(t,s,\hat{Z}(t,s),\hat{K}(t,s,\cdot))ds\\
\text{ \ \ \ \ \ \ \ \ }-\int_{t}^{T}\hat{Z}(t,s)dB(s)-\int_{t}^{T}\int_{%
\mathbb{R}
_{0}}\hat{K}(t,s,\zeta)\tilde{N}(ds,d\zeta),t\in\left[  0,T\right]  .
\end{array}
\]
Define%
\[%
\begin{array}
[c]{ll}%
\tilde{Y}(t) & :=\lambda Y^{-\psi_{1}(\cdot)}(t)+(1-\lambda)Y^{-\psi_{2}%
(\cdot)}(t),\\
\tilde{Z}(t,s) & :=\lambda Z^{-\psi_{1}(\cdot)}(t,s)+(1-\lambda)Z^{-\psi
_{2}(\cdot)}(t,s),\\
\tilde{K}(t,s,\cdot) & :=\lambda K^{-\psi_{1}(\cdot)}(t,s,\cdot)+(1-\lambda
)K^{-\psi_{2}(\cdot)}(t,s,\cdot).
\end{array}
\]
Then%
\[%
\begin{array}
[c]{l}%
\tilde{Y}(t)=-\lambda\psi_{1}(t)-(1-\lambda)\psi_{2}(t)\\
+\int_{t}^{T}[\lambda g(t,s,Z^{-\psi_{1}(\cdot)}(t,s),K^{-\psi_{1}(\cdot
)}(t,s,\cdot))+(1-\lambda)g(t,s,Z^{-\psi_{2}(\cdot)}(t,s),K^{-\psi_{2}(\cdot
)}(t,s,\cdot))]ds\\
\text{ \ \ \ \ \ \ \ \ }-\int_{t}^{T}\tilde{Z}(t,s)dB(s)-\int_{t}^{T}\int_{%
\mathbb{R}
_{0}}\tilde{K}(t,s,\zeta)\tilde{N}(ds,d\zeta)\\
\geq-\lambda\psi_{1}(t)-(1-\lambda)\psi_{2}(t)+\int_{t}^{T}g(t,s,\tilde
{Z}(t,s),\tilde{K}(t,s,\cdot))\\
-\int_{t}^{T}\tilde{Z}(t,s)dB(s)-\int_{t}^{T}\int_{%
\mathbb{R}
_{0}}\tilde{K}(t,s,\zeta)\tilde{N}(ds,d\zeta),t\in\left[  0,T\right]  ,
\end{array}
\]
where we have used the convexity of $g$ in the last inequality. By the
comparison Theorem \ref{comp}, we deduce that
\[
\hat{Y}(t)\leq\tilde{Y}(t),\text{for each }t\in\left[  0,T\right]  .
\]
Hence%
\begin{align*}
\rho(t;\lambda\psi_{1}(\cdot)+(1-\lambda)\psi_{2}(\cdot))  &  =\hat{Y}%
(t)\leq\tilde{Y}(t)\\
&  =\lambda Y^{-\psi_{1}(\cdot)}(t)+(1-\lambda)Y^{-\psi_{2}(\cdot)}(t)\\
&  =\lambda\rho(t;\psi_{1}(\cdot))+(1-\lambda)\rho(t;\psi_{2}(\cdot)).
\end{align*}
(ii) Monotonicity: If $\psi_{1}(\cdot)\leq\psi_{2}(\cdot)$, then, by the
comparison Theorem \ref{comp}, $Y^{-\psi_{2}(\cdot)}(t)\leq Y^{-\psi_{1}%
(\cdot)}(t).$ Consequently,
\[
\rho(t;\psi_{2}(\cdot))=Y^{-\psi_{2}(\cdot)}(t)\leq Y^{-\psi_{1}(\cdot
)}(t)=\rho(t;\psi_{1}(\cdot)).
\]
(iii) Translation invariant: If $\psi(\cdot)\in L_{\mathcal{F}_{T}}^{2}(0,T)$
and $a\in%
\mathbb{R}
$ is a real constant. Then we get that
\[
Y^{-\psi(\cdot)+a}(t)=Y^{-\psi(\cdot)}(t)+a.
\]
Thus,%
\begin{align*}
\rho(t;\psi(\cdot)+a)  &  =Y^{-\psi(\cdot)+a}(t)=Y^{-\psi(\cdot)}(t)-a\\
&  =\rho(t;\psi(\cdot))-a,\text{ for each }t\in\left[  0,T\right]  .
\end{align*}
(iv) The past independence is a direct consequence of the definition of $\rho
$. \hfill$\square$ \bigskip

\section{Semimartingale issues}

In this section, we will discuss some particular cases where the solution $Y$
of the above BSVIE can be a semimartingale. \newline For simplicity, we do not
consider jumps, since the jump terms do not play an essential role
here.\newline\newline Consider the solution couple $\left(  Y,Z\right)  \in
L_{y}^{2}\times L_{z}^{2}$ of a BSVIE of the form%
\begin{equation}%
\begin{array}
[c]{c}%
Y\left(  t\right)  =\psi\left(  t\right)  +%
{\textstyle\int_{t}^{T}}
g\left(  t,s,Y\left(  s\right)  ,Z\left(  t,s\right)  \right)  ds-%
{\textstyle\int_{t}^{T}}
Z\left(  t,s\right)  dB(s),\text{ }0\leq t\leq T,
\end{array}
\label{bsvie}%
\end{equation}
where $g:\triangle\times\mathbb{R}\times\mathbb{R}\times\Omega\rightarrow
\mathbb{R}$ is a Lipschitz driver and the terminal value $\psi\left(
t\right)  \in L_{\mathcal{F}_{T}}^{2}[0,T].$\newline In what follows, we
denote by the semimartingale $X(t)$ the solution of the stochastic
differential equation%
\begin{equation}
X(t)=x_{0}+%
{\textstyle\int_{0}^{t}}
b\left(  s,X(s)\right)  ds+%
{\textstyle\int_{0}^{t}}
\sigma\left(  s,X(s)\right)  dB(s),\text{ }t\in\lbrack0,T]. \label{sde}%
\end{equation}

\begin{center}
\textbf{Type 1 - BSVIE}
\end{center}

\noindent Let the couple $\left(  Y,Z\right)  \in L_{y}^{2}\times L_{z}^{2}$
be solution of the following BSVIE
\begin{equation}
Y(t)=F(X(t),X(T))-%
{\textstyle\int_{t}^{T}}
Z(t,s)dB(s), \label{y_1}%
\end{equation}
for some function $F:%
\mathbb{R}
^{2}\rightarrow%
\mathbb{R}
$ and $X(t)$ as given above by (\ref{sde}).\newline Now define%
\[
F(X(t),X(T)):=F_{1}(X(t))F_{2}(X(T)),
\]
for functions $F_{1}:%
\mathbb{R}
\rightarrow%
\mathbb{R}
$ and $F_{2}:%
\mathbb{R}
\rightarrow%
\mathbb{R}
$ which are assumed to be twice continuously differentiable ($C^{2}$).
Consider
\[
Y(t):=F_{1}(X(t))\tilde{Y}(t),
\]
where $\tilde{Y}(t)$ is the solution of the BSDE
\[
\tilde{Y}(t):=F_{2}(X(T))-%
{\textstyle\int_{t}^{T}}
Z(s)dB(s),
\]
and
\[
Z(t,s):=F_{1}(X(t))Z(s).
\]
By the Itô formula, we get that $Y(t)$ solution of (\ref{y_1}) is
a\ semimartingale.\newline

\begin{center}
\textbf{Type 2 - BSVIE}
\end{center}

\noindent Similarly as in the previous case, we consider a BSVIE of the form%
\begin{equation}
Y(t)=F(X(t),X(T))-%
{\textstyle\int_{t}^{T}}
Z(t,s)dB(s), \label{y}%
\end{equation}
\newline for functions $F\in C^{2}(%
\mathbb{R}
^{2})$. Then, for
\[
\tilde{Y}(t,x)=F(x,X(T))-%
{\textstyle\int_{t}^{T}}
\tilde{Z}(s,x)dB(s),
\]
we have that
\[
Y(t):=\tilde{Y}(t,X(t)),Z(t,s):=\tilde{Z}(s,X(t)).
\]
Using the Itô-Ventzell formula, we obtain that $Y(t)$ given by (\ref{y}) is a
semimartingale.\newline\newline

\begin{center}
\textbf{Type 3 - BSVIE}
\end{center}

\noindent Now we consider a BSVIE for a driver $g$ which does not depend on
$Z$, as follows:
\[
Y(t)=F(X(t),X(T))+%
{\textstyle\int_{t}^{T}}
g(X(t),X(s),Y(s))ds-%
{\textstyle\int_{t}^{T}}
Z(t,s)dB(s),t\in\lbrack0,T].
\]
Knowing $Y$, we can consider
\[
\tilde{Y}(t,x)=F(x,X(T))+%
{\textstyle\int_{t}^{T}}
g(x,X(s),Y(s))ds-%
{\textstyle\int_{t}^{T}}
\tilde{Z}(s,x)dB(s),t\in\lbrack0,T].
\]
Define
\begin{align*}
\bar{Y}(t)  &  :=\tilde{Y}(t,X(t)),\\
\bar{Z}(t,s)  &  :=\tilde{Z}(s,X(t)).
\end{align*}
Then
\[
\bar{Y}(t)=F(X(t),X(T))+%
{\textstyle\int_{t}^{T}}
g(X(t),X(s),Y(s))ds-%
{\textstyle\int_{t}^{T}}
\bar{Z}(t,s)dB(s),t\in\lbrack0,T].
\]
By uniqueness of the solution, we have
\begin{align*}
Y(t)  &  =\bar{Y}(t)=\tilde{Y}(t,X(t)),\\
Z(t,s)  &  =\bar{Z}(t,s)=\tilde{Z}(s,X(t)).
\end{align*}
By Itô-Ventzel's formula, we get that $Y(t)=\tilde{Y}(t,X(t))$ is a
semimartingale.\newline\newline The most general case, i.e., when the driver
depends on both $Y$ and $Z$ is still open and it is a subject of further research.

\bigskip

\noindent\textbf{Acknowledgement. }We would like to thank Prof. Yong for
pointing out the reference \cite{wy} which helped us to improve the paper. We
also want to thank Prof. Rosazza Gianin for helpful comments.

\end{document}